\documentclass[12pt]{article}
\usepackage[pdftex]{graphicx}

\usepackage{url}

\usepackage[unicode=true,bookmarks=true,bookmarksnumbered=false,
bookmarksopen=true,bookmarksopenlevel=1,breaklinks=false,
pdfborder={0 0 1},backref=section,colorlinks=false]{hyperref}

\usepackage{latexsym,amssymb,amsmath}
\usepackage{amsthm}

\renewcommand{\qed}{\hfill\rule[-0pt]{5pt}{15pt}}

\font\defi=cmr7
\newcommand{\adef}{\,\,\textrm{\raisebox{1.5pt}{\defi :}=}\,}

\newcommand{\aeq}{\,\,\textrm{\phantom{\raisebox{1.5pt}{\defi :}}=}\,}

\usepackage{amsopn}
\DeclareMathOperator{\re}{Re}

\newtheorem{thm}{Theorem}
\newtheorem{lemma}{Lemma}

\newcommand{\nn}{\mathbb N}

\newcommand{\rr}{\mathbb R}
\newcommand{\cc}{\mathbb C}

\newcommand{\modulo}[1]{\left\vert{#1}\right\vert}
\newcommand{\norma}[1]{\left\Vert{#1}\right\Vert}
\newcommand{\pent}[1]{\left\lfloor{#1}\right\rfloor}
\newcommand{\pfrac}[1]{\left\lbrace{#1}\right\rbrace}

\newcommand{\bin}[2]{\pfrac{\frac{x}{{#1}}}\;-\;\frac{{#2}}{{#1}}\,\pfrac{\frac{x}{{#2}}}}

\begin{document}

\DeclareGraphicsExtensions{.jpg,.pdf,.mps,.png,.svg}


\begin{center}

{\bf\Large An Algorithm to Generate Square-Free Numbers and to Compute 
  the M\"obius Function}

\vspace{5mm}

F. Auil\\{\tt auil@usp.br}

\end{center}

\vspace{0.2cm}


\begin{quote}
{\small
{\bf Abstract.}
We introduce an algorithm that iteratively produces a sequence of
natural numbers $k_{i}$ and functions $b_{i}$ defined in the interval
$[1,+\infty)$. The number $k_{i+1}$ arises as the first point of
discontinuity of $b_{i}$ above $k_{i}$.
We derive a set of properties of both sequences, suggesting
that (1) the algorithm produces square-free
numbers $k_{i}$, (2) \emph{all} the square-free numbers are generated
as the output of the algorithm, and (3) the value of the M\"obius
function $\mu(k_{i})$ can be evaluated as $b_{i}(k_{i+1}) -
b_{i}(k_{i})$.
The logical equivalence of these properties is rigorously
proved. The question remains open if one of these
properties can be derived from the definition of the algorithm.
Numerical evidence, limited to $5\times 10^{6}$, seems to support this 
conjecture.}

\bigskip
{\bf Keywords:} M\"obius function, square-free numbers, zeta function,
Riemann hypothesis.

\bigskip
{\bf AMS Subject Classification:} 11Y55, 11M99, 11Y35.

\end{quote}


\section{Introduction}

A natural number $n\in\nn$ is called {\it square-free}, if the
exponents arising in its prime factorization
\begin{equation*}
n = p_{1}^{r_{1}}\,p_{2}^{r_{2}}\cdots p_{k}^{r_{k}}
\end{equation*}
are all equal to 1, i.e., $r_{1}=r_{2}=\cdots=r_{k}=1$.
For a natural number $n$ with prime factorization as above, the {\it
  M\"obius function} $\mu$ is defined as
\begin{equation*}
\mu(n) =
\begin{cases}
1 & \text{if $n=1$;} \\
(-1)^{\sum_{i=1}^{k} r_{i}}, & \text{if $n$ is square-free;} \\
0, & \text{otherwise.}
\end{cases}
\end{equation*}
In other words, $\mu(n)$ is zero when $n$ has a square factor, and
otherwise gives the parity of the number of (distinct) prime factors
of $n$. The M\"obius function has important applications in number
theory, many of them concerning to the Riemann hypothesis about the
zeros of the zeta function, \cite{apo}, \cite{edw}, \cite{har},
\cite{tit}. \\

An outstanding problem in algorithmic number theory
is to compute $\mu(n)$ efficiently without first factoring $n$.  By
``efficiently'' we mean a number of bit operations bounded by a
polynomial in $\log n$, the length of $n$ in binary. As far as we
know, the question remains open if the computation of $\mu(n)$ can be
done in polynomial time, and in fact, nobody currently knows a way to
compute it significantly faster than factoring $n$. \\

In this paper, we present an algorithm that iteratively produces a
sequence of numbers $k_{i}$ and the value of $\mu(k_{i})$. 
In order to determine $\mu(k_{i})$, it is necessary to generate the
whole sequence $k_{1}, k_{2},\dots, k_{i}$. Our algorithm is based on
a sequence of arithmetical functions $b_{i}$, with the numbers $k_{i}$
arising as discontinuity points. These functions are closely related
to the Nyman-Beurling approach to the Riemann hypothesis. The main
result of the present paper is a set of properties of the sequences
$k_{i}$ and $b_{i}$, suggesting that
\begin{itemize}
\item
The numbers $k_{i}$ generated by the algorithm are square-free.
\item
The set of \emph{all} the square-free numbers can be generated by the
algorithm.
\item
The value of the M\"obius function $\mu(k_{i})$ can be evaluated as
$b_{i}(k_{i+1}) - b_{i}(k_{i})$.
\end{itemize}
We are able to prove the logical equivalence of these
properties. Unfortunately, using the definition of the algorithm, we
cannot prove, neither disprove, if one (then all) of these conditions
are satisfied. Numerical evidence, limited to $5\times 10^{6}$, seems
to support the conjectures quoted above. \\


This paper is organized as follows. The section \ref{hilbert}
introduces the framework of the Nyman-Beurling approach to Riemann
hypothesis. The section \ref{basic} provides the basic definitions for
this paper, and the algorithm is defined in section
\ref{algorithm}. The main result of the present work is stated and
proved in section \ref{main}, while further remarks are done in
section \ref{discussion}. Finally, the implementation of the algorithm
and numerical results are discussed in section \ref{numerical}.


\section{Hilbert Space Approach to Riemann Hypothesis}
\label{hilbert}

Denote by $\pent{x}$ the integer part of $x$, i.e., the greatest
integer less than, or equal to, $x$. Define the {\it fractional part}
function by $\pfrac{x}=x - \pent{x}$. 
Given $n\in\nn$ and two families of parameters
$\{a_{k}\}_{k=1}^{n}\subset\cc$ and $\{\theta_{k}\}_{k=1}^{n}\subset
(0,1]$, we define a \emph{Beurling function} as a function
$F_{n}$ (the sub-index $n$ included in the notation for
convenience) of the form
\begin{equation}
F_{n}(x)\adef\sum_{k=1}^{n} a_{k}\,\pfrac{\frac{\theta_{k}}{x}}.
\label{eq-1}
\end{equation} 

\noindent
For a Beurling function $F_{n}$, an elementary computation shows that 
\begin{equation}
\int_{0}^{1} (F_{n}(x)+1) \, x^{s-1} \, dx 
= \frac{\displaystyle\sum_{k=1}^{n} a_{k}\,\theta_{k}}{s-1}
+ \frac{1}{s}\left(1-\zeta(s) \, \sum_{k=1}^{n}
  a_{k}\,\theta_{k}^{s}\right);
\label{eq-2}
\end{equation}

\noindent
for the complex variable $s\in\cc$ in the half-plane $\re(s)>0$, i.e.
for $s$ with positive real part. Here, $\zeta$ denotes the Riemann's
{\it zeta function} given by
\begin{equation*}
\zeta(s)\adef\sum_{k=1}^{\infty}\frac{1}{k^{s}}.
\end{equation*}
The classic references are \cite{tit} and \cite{edw}.
A derivation of the relation (\ref{eq-2}) can be found, for instance,
in \cite[p. 253]{don}. It is useful (but
not always necessary from a theoretical point of view) to assume that
the parameters defining the function $F_{n}$ satisfy the additional
condition 
\begin{equation}
\sum_{k=1}^{n} a_{k}\,\theta_{k}=0. 
\label{eq-3}
\end{equation}

\noindent
In this case, the first term at the right-hand side of (\ref{eq-2})
vanishes, simplifying the expression. The identity (\ref{eq-2}) is the
starting point of the following theorem by Beurling


\begin{thm}[Beurling]

The zeta function $\zeta(s)$ has no zeros in the half-plane $\re(s) >
1/p$ if and only if the set of (Beurling) functions $\{f_{\theta}(x) =
\pfrac{\theta/x}\}_{0<\theta\leqslant 1}$ is dense in
$L^{p}([0,1],dx)$. 

\end{thm}


\noindent
See \cite[p. 252]{don} for a proof of the Beurling theorem and further
references. Note that for $p=2$ this result provides an equivalent
condition to the {\it Riemann hypothesis} (RH) for the zeta
function. This is the Beurling, or Nyman-Beurling, approach to RH. \\

The Beurling theorem above has an easy half part, whose proof can be
sketched as follows. From relation (\ref{eq-2}) and assuming
(\ref{eq-3}), we have
\begin{equation*}
\modulo{\frac{1}{s}\left(1-\zeta(s)\,\sum_{k=1}^{n}a_{k}\,\theta_{k}^{s}\right)}
=\modulo{\int_{0}^{1} (F_{n}(x)+1)\,x^{s-1}\,dx}
\leqslant\norma{F_{n}(x)+1}\,\norma{x^{s-1}};
\end{equation*}
where the last relation follows using the Schwarz inequality in 
$L^{2}([0,1],dx)$. Therefore, if the first norm in the
right-hand side above can be done arbitrarily small for a suitable
choice of $n$, $a_{k}$'s and $\theta_{k}$'s, then the function
$\zeta(s)$ could not have zeros for $\re(s) > 1/2$. We will refer
to the first condition above as the {\it Beurling criterion} (BC) for
RH.  Note also that in order to demonstrate the RH, it is sufficient to
prove that the constant function equal to $-1$ can be arbitrarily
approximated in the norm of the Hilbert space $L^{2}([0,1],dx)$ by
Beurling functions $F_{n}$ of the form (\ref{eq-1}). It was proved in
\cite{bae2} that BC remains equivalent to RH if the parameters
$\theta_{k}$ are restricted to be reciprocal of natural numbers,
i.e. $\theta_{k}=1/b_{k}$, with $b_{k}\in\nn$.\\


Several approximating functions to $-1$ of the form (\ref{eq-1}) were
proposed in the literature. From the relations (\ref{eq-2}) and
(\ref{eq-3}), we have that under the BC, the ``partial sum''
\begin{equation}
\sum_{k=1}^{n} a_{k}\,\theta_{k}^{s}.
\end{equation}
is an approximation to the reciprocal of the zeta function 
$1/\zeta(s)$, which is known to have an expression as a Dirichlet
series  
\begin{equation}
\frac{1}{\zeta(s)}\;=\;\sum_{k=1}^{n} \frac{\mu(k)}{k^{s}}, 
\end{equation}
convergent for $\re(s) > 1$. Therefore, a (naive) first choice for an
approximating function would be
\begin{equation}
S_{n}(x)\adef\sum_{k=1}^{n}\mu(k)\,\pfrac{\frac{1/k}{x}}.
\label{eq-6}
\end{equation}
Note that this function does not matches the condition
(\ref{eq-3}). We can handle this without subtlety, just by
subtracting the difference, that is given by $g(n)$, where
\begin{equation}
g(t)\adef \sum_{\nn\ni k\leqslant t}\frac{\mu(k)}{k}. 
\label{eq-6b}
\end{equation}
Therefore, a second choice would be
\begin{align}
B_{n}(x) & \adef\sum_{k=1}^{n} \mu(k)\,\pfrac{\frac{1/k}{x}} \;
- \;n\,g(n)\,\pfrac{\frac{1/n}{x}} \\
 & \aeq\sum_{k=1}^{n-1} \mu(k)\,\pfrac{\frac{1/k}{x}} \;
- \;n\,g(n-1)\,\pfrac{\frac{1/n}{x}}. 
\label{eq-7}
\end{align}
Other variants were also proposed, as
\begin{equation}
V_{n}(x)\adef\sum_{k=1}^{n} \mu(k)\,\pfrac{\frac{1/k}{x}} \;
- \; g(n)\,\pfrac{\frac{1}{x}}.
\label{eq-8}
\end{equation}
Unfortunately, the sequences (\ref{eq-6}), (\ref{eq-7}) and
(\ref{eq-8}) are known to be \emph{not} convergent to $-1$ in
$L^{2}([0,1],dx)$, as proved in \cite{bae}. A survey on the
Nyman-Beurling reformulation of the Riemann hypothesis
and later developments by Baez-Duarte can be found in \cite{bag}.



\section{Basic Definitions}
\label{basic}

In order to motivate the definitions below, assume that a Beurling
function $F_{n}$ as in (\ref{eq-1}) is constant between the reciprocal
of the natural numbers. In other words, assume
that such a function takes a constant value in each of the intervals
$(\frac{1}{k+1},\frac{1}{k}]$, for all $k\in\nn$ (but the constant
value may differ from interval to interval).
In this case, the integral in (\ref{eq-2}) can be expressed
alternatively as an infinite series, involving the values of
$f_{n}(k) \adef F_{n}(1/k)$. Furthermore, the values of $F_{n}(x)$ for
all $x\in[0,1]$ are completely determined by the values of
$f_{n}(k)$ for all $k\in\nn$. We will call an {\it arithmetical
 Beurling function} a function of the form $f_{n}(k) = F_{n}(1/k)$,
where $F_{n}$ is a Beurling function. \\

We introduce now, perhaps the simplest, non-trivial, example of
arithmetical Beurling function satisfying the condition (\ref{eq-3}).
For $a,b\in\nn$ define the function $\beta_{a,b}$ as
\begin{equation}
\beta_{a,b}(x)\adef\bin{a}{b}, 
\end{equation}
As the functions $\beta_{a,b}$ will be the basic blocks
in our construction, we summarize some of its elementary properties in
the following result. 


\begin{lemma}
Consider $a,b\in\rr$, with $0<a<b$. Then,

\begin{itemize}

\item[\bf{a.}]
$\displaystyle{\pfrac{\frac{x}{a}}}$ and 
$\displaystyle{\pfrac{\frac{x}{b}}}$ are right-continuous, and
linearly independent functions.

\item[\bf{b.}]
$\beta_{a,b}(x)=0$, when $0\leqslant x<a$. 

\item[\bf{c.}]
Let $k\in\nn$ be such that $(k-1)a<b\leqslant ka$. Then,  
$$
\beta_{a,b}(x)=\begin{cases}-j & \text{\rm if $ja\leqslant x<(j+1)a$, for
$j=1,\dots,(k-2)$;}\\-(k-1) & \text{\rm if $(k-1)a\leqslant x<b$.}\end{cases}
$$

\item[\bf{d.}]
Assume $a,b\in\nn$. Then, $\beta_{a,b}(x)$ is constant when
$k\leqslant x<k+1$, for all $k\in\nn$. 

\end{itemize}
\label{lemma-1}
\end{lemma}


\paragraph{\it Proof:}
(a): The right-continuity is derived from of $\pent{x}$. Now, if 
$c_{1}\,\pfrac{\frac{x}{a}}\;+\;c_{2}\,\pfrac{\frac{x}{b}}=0$ for all
$x$, then for $x=a$ we have  
$0=c_{1}\,\pfrac{\frac{a}{a}}\;+\;c_{2}\,\pfrac{\frac{a}{b}}
=c_{2}\,\frac{a}{b}$.
Thus, $c_{2}=0$ and we have $c_{1}\,\pfrac{\frac{x}{a}}=0$ for all
$x$, and taking now $x=a/2$ we get 
$0=c_{1}\,\pfrac{\frac{1}{2}}=c_{1}/2$ and $c_{1}=0$. \\


\noindent
(b): If $0\leqslant x<a<b$, then $x/b<1$ and $x/a<1$. Thus, 
$\pfrac{\frac{x}{a}} - \frac{b}{a}\pfrac{\frac{x}{b}}
= \frac{x}{a} - \frac{b}{a}\frac{x}{b} = 0$.\\


\noindent
(c): Assume $j=1,\dots,(k-2)$. Then, for $ja\leqslant x < (j+1)a <b$,
we have $x/b < 1$ and $j \leqslant x/a < (j+1)$. Thus, 
$\pfrac{\frac{x}{a}} - \frac{b}{a}\pfrac{\frac{x}{b}}
= \frac{x}{a} - \pent{\frac{x}{a}}-\frac{b}{a}\frac{x}{b}
= \frac{x}{a} - j - \frac{x}{a}
= -j$. 
Analogously, for $(k-1)a \leqslant x < b < ka$,
we have $x/b < 1$ and $(k-1) \leqslant x/a < k$. Thus, 
$\pfrac{\frac{x}{a}} - \frac{b}{a}\pfrac{\frac{x}{b}}
= \frac{x}{a} - \pent{\frac{x}{a}} - \frac{b}{a}\frac{x}{b}
= \frac{x}{a} - (k-1) - \frac{x}{a}
= -(k-1)$. \\


\noindent
(d): If $x<b$ then (d) is true by (b) and (c) already
proven. Consider now  $b \leqslant x$. If $k < x < k+1$ then
$x\notin\nn$ and thus $x/a\notin\nn$ and  $x/b\notin\nn$. Therefore,
there exists $k_{0}$ and $l_{0}$ in $\nn$ such that 
$ak_{0} < x < a(k_{0}+1)$ and $bl_{0} < x <b(l_{0}+1)$, and we have 
$\pfrac{\frac{x}{a}} - \frac{b}{a}\pfrac{\frac{x}{b}}
= \frac{x}{a} - k_{0} - \frac{b}{a}(\frac{x}{b}-l_{0})
= -k_{0} + \frac{b}{a}l_{0}$, which is a constant independent of $x$. 
As $\pfrac{\frac{x}{a}} - \frac{b}{a}\pfrac{\frac{x}{b}}$ is
right-continuous, by part (a), this is also true for 
$k \leqslant x < k+1$. \qed


\section{The Algorithm}
\label{algorithm}

We will define a sequence of numbers $\{k_{i}\}_{i\in\nn}$ and functions
$\{b_{i}\}_{i\in\nn}$ iteratively as follows. Start with the following
definitions 
\begin{equation}
\begin{split}
k_{1} & \adef 1;\\
k_{2} & \adef 2;\\
b_{2}(x) & \adef\pfrac{\frac{x}{k_{1}}}\; 
- \; \frac{k_{2}}{k_{1}}\,\pfrac{\frac{x}{k_{2}}}.
\label{eq-20}
\end{split}
\end{equation}
Assuming now that $k_{i}$ and $b_{i}$ are already defined,
for $i\geqslant 2$ define iteratively $k_{i+1}$ and $b_{i+1}$ as
follows. The number $k_{i+1}$ is defined as $k_{i+1}\adef k_{i}+j$,
where $j$ is the least integer such that $b_{i}(k_{i}+j)\neq
b_{i}(k_{i})$. Once determined the number $k_{i+1}$, the function
$b_{i+1}$ is defined as 
\begin{equation}
b_{i+1}(x)
\adef b_{i}(x) \;
+ \; \left(1 + b_{i}(k_{i})\right)\left(\pfrac{\frac{x}{k_{i}}}\;
- \; \frac{k_{i+1}}{k_{i}}\,\pfrac{\frac{x}{k_{i+1}}}\right).
\label{eq-21}
\end{equation}
Some elementary properties derived from these definitions are
summarized in the following result.


\begin{lemma}

For any $i\in\nn$ we have

\begin{itemize}

\item[{\bf a.}] $b_{i}$ is a right-continuous function, which is
  constant between the natural numbers.

\item[{\bf b.}] $b_{i+1}(k_{i})=-1$.

\item[{\bf c.}] Assume $k_{i+1}\leqslant 2k_{i}$ for $i\geqslant
  2$. Then, $b_{i}(x) = -1$ for all $x \in [1,k_{i})$. In particular,
  the sequence $\{b_{i}\}_{i\in\nn}$ converges point-wise to $-1$ in
$[1,+\infty)$. 

\end{itemize}
\label{lemma-2}
\end{lemma}


\paragraph{\it Proof:}
(a): Observe that each $b_{i}$ is a (finite) linear combination of
$\beta_{p,q}$. Therefore, this result is a direct consequence of Lemma
\ref{lemma-1}. \\


\noindent
(b): Is an immediate consequence of definition (\ref{eq-21}). \\


\noindent
(c): For induction on $i$. The case $i=2$ is an immediate
consequence of definition (\ref{eq-20}). Assume now that $b_{j}(x) =
-1$ when $x \in [1,k_{j})$ for all $j \leqslant i$. If $x < k_{i} <
k_{i+1}$, then from definition (\ref{eq-21}) we have $b_{i+1}(x) =
b_{i}(x)$ which is equal to $-1$ by the inductive hypothesis. Now if 
$k_{i} \leqslant x < k_{i+1} \leqslant 2k_{i}$,  also from definition
(\ref{eq-21}) we have
\begin{equation*}
b_{i+1}(x)
= b_{i}(x) 
+ (1+b_{i}(k_{i}))\,\left(\frac{x}{k_{i}} - 1 - \frac{k_{i+1}}{k_{i}}\,\frac{x}{k_{i+1}}\right)
= b_{i}(x) - 1 - b_{i}(k_{i})
= -1,
\end{equation*}
because by the definition of $k_{i+1}$, the function $b_{i}(x)$ is
constant for $x\in[k_{i},k_{i+1})$. \qed


\paragraph{Remark:} As $b_{i}$ is a constant function between the
natural numbers, the number $k_{i+1}$ is the first point of
discontinuity of $b_{i}$ above $k_{i}$.


\section{Main Result}
\label{main}

The next result is relevant in order to establish a relationship between
the sequence $\{k_{i}\}_{i\in\nn}$, the values of $\mu(k_{i})$, and
the square-free numbers.


\begin{lemma}
The following conditions are equivalent

\begin{itemize}

\item[\bf{a.}]
$\displaystyle{
\sum_{j=1}^{i}\frac{\mu(k_{j})}{k_{j}}\;=\;\frac{1+b_{i}(k_{i})}{k_{i}}
}$, for $i\geqslant 2$.

\item[\bf{b.}]
$\displaystyle{
b_{i}(x)=\sum_{j=1}^{i-1}\mu(k_{j})\,\pfrac{\frac{x}{k_{j}}} \; 
- \;
k_{i}\left(\sum_{j=1}^{i-1}\frac{\mu(k_{j})}{k_{j}}\right)\pfrac{\frac{x}{k_{i}}}
}$,
for $i\geqslant 2$.

\item[\bf{c.}]
$\displaystyle{
\frac{\mu(k_{i})}{k_{i}}=\frac{1+b_{i}(k_{i})}{k_{i}}\;
- \; \frac{1+b_{i-1}(k_{i-1})}{k_{i-1}}
}$, for $i\geqslant 3$.

\end{itemize}

\noindent
Furthermore, if the condition $k_{i+1} < 2k_{i}$ is valid for
$i\geqslant 2$, then all conditions above are also equivalent to the
following ones

\begin{itemize}

\item[\bf{d.}]
$\displaystyle{
\mu(k_{i+1})=b_{i}(k_{i+1}) \; - \; b_{i}(k_{i})
}$, for $i\geqslant 2$. 

\item[\bf{e.}]
$\displaystyle{
\sum_{j=1}^{i}\mu(k_{j})\pent{\frac{k_{i}}{k_{j}}}=1
}$, for $i\geqslant 1$.

\end{itemize}
\label{lemma-3}
\end{lemma}


\paragraph{\it Proof:} In order to prove the logical equivalence
between all conditions in Lemma \ref{lemma-3}, we separately will
prove, first of all, the equivalences (a)$\Leftrightarrow$(b) and
(a)$\Leftrightarrow$(c). Then, after the introduction of the
additional condition $k_{i+1}<2k_{i}$, we will prove the equivalence
between (d)$\Leftrightarrow$(c) and (e)$\Leftrightarrow$(b).\\


\noindent
(a)$\Rightarrow$(b): For induction on $i$. For $i=2$, from
definition (\ref{eq-20}) we have 
\begin{equation}
b_{2}(x)
= \pfrac{\frac{x}{k_{1}}} - \frac{k_{2}}{k_{1}}\,\pfrac{\frac{x}{k_{2}}}
= \mu(k_{1})\,\pfrac{\frac{x}{k_{1}}} -
k_{2}\,\frac{\mu(k_{1})}{k_{1}}\,\pfrac{\frac{x}{k_{2}}}.
\end{equation}
Assuming now that condition (b) follows for all $j$ such that
$j\leqslant i$, we have
\begin{multline}
b_{i+1}(x)
= b_{i}(x) 
+ \left(1+b_{i}(k_{i})\right)
  \left(\pfrac{\frac{x}{k_{i}}} - \frac{k_{i+1}}{k_{i}}\,\pfrac{\frac{x}{k_{i+1}}}\right)\\
= \sum_{j=1}^{i-1}\mu(k_{j})\,\pfrac{\frac{x}{k_{j}}}
- k_{i}\left(\sum_{j=1}^{i-1}\frac{\mu(k_{j})}{k_{j}}\right)\pfrac{\frac{x}{k_{i}}}
+ \left(1+b_{i}(k_{i})\right)
\left(\pfrac{\frac{x}{k_{i}}} - \frac{k_{i+1}}{k_{i}}\,\pfrac{\frac{x}{k_{i+1}}}\right) \\
= \sum_{j=1}^{i-1}\mu(k_{j})\,\pfrac{\frac{x}{k_{j}}}
- k_{i}\left(\sum_{j=1}^{i-1}\frac{\mu(k_{j})}{k_{j}}\right)\pfrac{\frac{x}{k_{i}}}
+ \left(k_{i}\sum_{j=1}^{i}\frac{\mu(k_{j})}{k_{j}}\right)
\left(\pfrac{\frac{x}{k_{i}}} - \frac{k_{i+1}}{k_{i}}\,\pfrac{\frac{x}{k_{i+1}}}\right)\\
= \sum_{j=1}^{i-1}\mu(k_{j})\,\pfrac{\frac{x}{k_{j}}}
- k_{i}\left(\sum_{j=1}^{i-1}\frac{\mu(k_{j})}{k_{j}}\right)\pfrac{\frac{x}{k_{i}}}\hfill\\
\hfill + k_{i}\left(\sum_{j=1}^{i}\frac{\mu(k_{j})}{k_{j}}\right)\pfrac{\frac{x}{k_{i}}}
- k_{i+1}\left(\sum_{j=1}^{i}\frac{\mu(k_{j})}{k_{j}}\right)\pfrac{\frac{x}{k_{i+1}}}\\ 
= \sum_{j=1}^{i-1}\mu(k_{j})\,\pfrac{\frac{x}{k_{j}}}
+ k_{i}\,\frac{\mu(k_{i})}{k_{i}}\,\pfrac{\frac{x}{k_{i}}}
- k_{i+1}\left(\sum_{j=1}^{i}\frac{\mu(k_{j})}{k_{j}}\right)\pfrac{\frac{x}{k_{i+1}}}\\
= \sum_{j=1}^{i}\mu(k_{j})\,\pfrac{\frac{x}{k_{j}}}
- k_{i+1}\left(\sum_{j=1}^{i}\frac{\mu(k_{j})}{k_{j}}\right)\pfrac{\frac{x}{k_{i+1}}};
\end{multline}
and this proves condition (b) for $i+1$. 
Here, in the first equality we have used the inductive hypothesis and
in the second one we have used condition (a).\\


\noindent
(b)$\Rightarrow$(a): From Lemma \ref{lemma-2} (b), by condition (b)
we have
\begin{equation}
\begin{split}
0
& = 1  +b_{i+1}(k_{i}) \\
& = 1 + \sum_{j=1}^{i}\mu(k_{j})\,\pfrac{\frac{k_{i}}{k_{j}}}
- k_{i+1}\left(\sum_{j=1}^{i}\frac{\mu(k_{j})}{k_{j}}\right)\pfrac{\frac{k_{i}}{k_{i+1}}}\\
& = 1 + \sum_{j=1}^{i-1}\mu(k_{j})\,\pfrac{\frac{k_{i}}{k_{j}}}
+ \mu(k_{i})\,\pfrac{\frac{k_{i}}{k_{i}}}
- k_{i+1}\left(\sum_{j=1}^{i}\frac{\mu(k_{j})}{k_{j}}\right)\frac{k_{i}}{k_{i+1}}\\
& = 1 + \sum_{j=1}^{i-1}\mu(k_{j})\,\pfrac{\frac{k_{i}}{k_{j}}}
- k_{i}\left(\sum_{j=1}^{i}\frac{\mu(k_{j})}{k_{j}}\right) \\
& = 1 + b_{i}(k_{i}) - k_{i}\left(\sum_{j=1}^{i}\frac{\mu(k_{j})}{k_{j}}\right);
\end{split}
\end{equation}
and this proves condition (a). Here we have used that 
$\pfrac{\frac{k_{i}}{k_{i+1}}} = \frac{k_{i}}{k_{i+1}}$, 
(because $\frac{k_{i}}{k_{i+1}} < 1$), that  
$\pfrac{\frac{k_{i}}{k_{i}}} = 0$, and also the relation
$b_{i}(k_{i}) = \sum_{j=1}^{i-1}\mu(k_{j})\,\pfrac{\frac{k_{i}}{k_{j}}}$,
which is an easy consequence of condition (b).\\ 


\noindent
(a)$\Rightarrow$(c): Using condition (a) for $i$ and $i-1$ we have
\begin{align}
\sum_{j=1}^{i}\frac{\mu(k_{j})}{k_{j}} & =
\frac{1+b_{i}(k_{i})}{k_{i}},\label{eq-25}\\
\sum_{j=1}^{i-1}\frac{\mu(k_{j})}{k_{j}} & =
\frac{1+b_{i-1}(k_{i-1})}{k_{i-1}};
\label{eq-26}
\end{align}
and subtracting (\ref{eq-26}) from (\ref{eq-25}) we get condition
(c).\\


\noindent
(c)$\Rightarrow$(a): Denoting
$\alpha(i)\adef\frac{1+b_{i}(k_{i})}{k_{i}}$, from (c) we have
\begin{multline}
\sum_{j=1}^{i}\frac{\mu(k_{j})}{k_{j}} 
=\frac{\mu(k_{1})}{k_{1}} \; + \; \frac{\mu(k_{2})}{k_{2}} \; + \;\sum_{j=3}^{i}\alpha(j)-\alpha(j-1) \;
= \; 1 \; - \;\frac{1}{2} \; + \; \alpha(i)-\alpha(2) \\
= \; 1 \; - \; \frac{1}{2} \; + \; \alpha(i) \; -  \; \frac{1+b_{2}(k_{2})}{k_{2}} \;
= \; 1 \; - \; \frac{1}{2} \; + \; \alpha(i) \; - \; \frac{1}{2} \;
= \; \alpha(i) \;
= \; \frac{1+b_{i}(k_{i})}{k_{i}}. 
\end{multline}
Here we have used that $b_{2}(k_{2})=0$, by definition (\ref{eq-20}).\\

\noindent
Therefore, we have proved
(b)$\Leftrightarrow$(a)$\Leftrightarrow$(c). Assume now condition
$k_{i+1}<2k_{i}$, for $i\geqslant 2$.\\


\noindent
(d)$\Leftrightarrow$(c): From definition (\ref{eq-21}) we have
\begin{equation}
\begin{split}
b_{i}(k_{i})
& = b_{i-1}(k_{i}) 
+ \left(1+b_{i-1}(k_{i-1})\right)
  \left(\pfrac{\frac{k_{i}}{k_{i-1}}} 
   - \frac{k_{i}}{k_{i-1}}\,\pfrac{\frac{k_{i}}{k_{i}}}\right) \\
& = b_{i-1}(k_{i})
+ \left(1+b_{i-1}(k_{i-1})\right)\,\pfrac{\frac{k_{i}}{k_{i-1}}} \\
& = b_{i-1}(k_{i}) 
+ \left(1+b_{i-1}(k_{i-1})\right)
  \left(\frac{k_{i}}{k_{i-1}} - \pent{\frac{k_{i}}{k_{i-1}}}\right) \\
& = b_{i-1}(k_{i}) 
+ \left(1+b_{i-1}(k_{i-1})\right)\left(\frac{k_{i}}{k_{i-1}} - 1 \right) \\
& = b_{i-1}(k_{i}) 
+ \left(1+b_{i-1}(k_{i-1})\right)\,\frac{k_{i}}{k_{i-1}} 
- \left(1+b_{i-1}(k_{i-1})\right). 
\label{eq-27}
\end{split}
\end{equation}
Observe that the additional condition implies
$1< \; \frac{k_{i}}{k_{i-1}} < 2$, and therefore
$\pent{\frac{k_{i}}{k_{i-1}}} = 1$. 
From (\ref{eq-27}) follows 
\begin{equation}
\frac{1+b_{i}(k_{i})}{k_{i}} \; - \;
\frac{1+b_{i-1}(k_{i-1})}{k_{i-1}} \;
= \; \frac{b_{i-1}(k_{i})-b_{i-1}(k_{i-1})}{k_{i}},
\label{eq-28}
\end{equation}
for $i \leqslant 3$. The equivalence between (d) and (c) is a direct
consequence of (\ref{eq-28}) above. \\


\noindent
(b)$\Rightarrow$(e): From Lemma \ref{lemma-2} (b), by condition (b)
we have
\begin{multline}
0 
= 1 + b_{i}(k_{i-1}) 
= 1 + \sum_{j=1}^{i-1}\mu(k_{j})\,\pfrac{\frac{k_{i-1}}{k_{j}}} 
- k_{i}\left(\sum_{j=1}^{i-1}\frac{\mu(k_{j})}{k_{j}}\right)\pfrac{\frac{k_{i-1}}{k_{i}}}\\
= 1 
+ \left(\sum_{j=1}^{i-1}\frac{\mu(k_{j})}{k_{j}}\right)\,k_{i-1} 
- \sum_{j=1}^{i-1}\mu(k_{j})\,\pent{\frac{k_{i-1}}{k_{j}}}
-
k_{i}\left(\sum_{j=1}^{i-1}\frac{\mu(k_{j})}{k_{j}}\right)\frac{k_{i-1}}{k_{i}} \\
+
k_{i}\left(\sum_{j=1}^{i-1}\frac{\mu(k_{j})}{k_{j}}\right)\pent{\frac{k_{i-1}}{k_{i}}}
= 1 - \sum_{j=1}^{i-1}\mu(k_{j})\,\pent{\frac{k_{i-1}}{k_{j}}},  
\end{multline}
and this proves condition (e). 
Here, we have used that $\frac{k_{i-1}}{k_{i}} < 1$ and therefore 
$\pent{\frac{k_{i-1}}{k_{i}}} = 0$. \\


\noindent
(e)$\Rightarrow$(b): Assume condition (e) valid for all
$i\in\nn$. We will prove (b) by induction on $i$. Condition (b) for
$i=2$ follows from definition (\ref{eq-20}) as done in the proof
(a)$\Rightarrow$(b) above. Assume now condition (b) valid for all
$n\leqslant i$. Using the notation 
$g(i) = \sum_{j=1}^{i}\frac{\mu(k_{i})}{k_{i}}$, we have
\begin{equation}
\begin{split}
b_{i}(k_{i+1}) 
& = \sum_{j=1}^{i-1}\mu(k_{j})\,\pfrac{\frac{k_{i+1}}{k_{j}}}
- k_{i}\left(\sum_{j=1}^{i-1}\frac{\mu(k_{j})}{k_{j}}\right)\pfrac{\frac{k_{i+1}}{k_{i}}}\\
& = k_{i+1}\,g(i-1) 
- \sum_{j=1}^{i-1}\mu(k_{j})\pent{\frac{k_{i+1}}{k_{j}}} 
- k_{i}\,g(i-1)\,\frac{k_{i+1}}{k_{i}}
+ k_{i}\,g(i-1)\,\pent{\frac{k_{i+1}}{k_{i}}} \\
& =
- \sum_{j=1}^{i-1}\mu(k_{j})\pent{\frac{k_{i+1}}{k_{j}}} 
+ k_{i}\,g(i-1) \\
& =
- \left(1 - \mu(k_{i+1})\pent{\frac{k_{i+1}}{k_{i+1}}}
- \mu(k_{i})\pent{\frac{k_{i+1}}{k_{i}}}\right)
+ k_{i}\,g(i-1) \\
& =
- \left(1 - \mu(k_{i+1}) - \mu(k_{i}) \right)
+ k_{i}\,g(i-1) \\
& =  k_{i}\,g(i-1) + \mu(k_{i+1}) + \mu(k_{i}) - 1.
\label{eq-29}
\end{split}
\end{equation}
Here we have used $\pent{\frac{k_{i+1}}{k_{i}}} = 1$ (a
consequence of the additional condition) and (e). Analogously we have 
\begin{equation}
\begin{split}
b_{i}(k_{i}) 
& = \sum_{j=1}^{i-1}\mu(k_{j})\,\pfrac{\frac{k_{i}}{k_{j}}}
-
k_{i}\left(\sum_{j=1}^{i-1}\frac{\mu(k_{j})}{k_{j}}\right)\pfrac{\frac{k_{i}}{k_{i}}} \\
& = \sum_{j=1}^{i-1}\mu(k_{j})\,\pfrac{\frac{k_{i}}{k_{j}}} \\
& = k_{i}\,g(i-1) 
- \sum_{j=1}^{i-1}\mu(k_{j})\pent{\frac{k_{i}}{k_{j}}} \\
& = k_{i}\,g(i-1) 
- \left(1 - \mu(k_{i})\pent{\frac{k_{i}}{k_{i}}}\right) \\
& = k_{i}\,g(i-1) - (1 - \mu(k_{i})) \\
& = k_{i}\,g(i-1) + \mu(k_{i}) - 1.
\label{eq-29b}
\end{split}
\end{equation}
Now, subtracting (\ref{eq-29b}) from (\ref{eq-29}) we get condition
(d). But we have already proved that
(d)$\Rightarrow$(c)$\Rightarrow$(a)$\Rightarrow$(b).\qed


\section{Discussion}
\label{discussion}

The definitions in section \ref{algorithm} provide an algorithm to
produce iteratively a sequence of numbers $\{k_{i}\}_{i\in\nn}$ and
functions $\{b_{i}\}_{i\in\nn}$. The Lemma \ref{lemma-3} in section
\ref{main} states a circle of logically equivalent properties of both
sequences. \\

The condition (b) of Lemma \ref{lemma-3} suggests that the functions
$b_{i}$ are the arithmetical counterpart of the approximating
functions $B_{n}$ in the relation (\ref{eq-7}). However, note that our
definition of $b_{i}$ in section \ref{algorithm} is quite different of
$B_{n}$. \\

The condition (d) of Lemma \ref{lemma-3} is related with the
conjecture that the algorithm produces square-free numbers
$k_{i}$, and also provides the value of the M\"obius function
$\mu(k_{i})$. Note that the value of $b_{i}(k_{i+1}) - b_{i}(k_{i})$
is never zero, by the definition of the algorithm in section
\ref{algorithm}. \\

The condition (e) of Lemma \ref{lemma-3} is related with the
conjecture that the algorithm produces \emph{all} the square-free
numbers. Indeed, if condition (e) were true, this would be a corollary
of a well known result; see \cite[p. 66]{apo}. \\

The condition $k_{i+1}<2k_{i}$ in Lemma \ref{lemma-3}, sufficient for
(d) and (e), seems to be also necessary, as the following heuristic
argument suggests. It is known that the square-free numbers are
distributed in $\nn$ with density $6/\pi^{2}$; see \cite[Thm. 333,
p. 269]{har}. Therefore, we can estimate the average distance between
two consecutive square-free numbers as $\pi^{2}/6$. Consequently,
$k_{i+1}\approx k_{i}+\pi^{2}/6$, or $k_{i+1}/k_{i}\approx
1+\pi^{2}/6k_{i}$. The last expression is less than $2$ for
$k_{i}>\pi^{2}/6\approx 1.64$. Thus, condition $k_{i+1}<2k_{i}$ for
$i\geqslant 2$ seems to be reasonable also. \\

Unfortunately, using the definition of the algorithm given in section
\ref{algorithm}, we cannot prove, neither disprove, if one of the
conditions in Lemma \ref{lemma-3} are satisfied. \\

Note also that the algorithm cannot compute isolated values of
$\mu(k_{i})$. In order to determine $\mu(k_{i})$, it is necessary to
generate the whole sequence $k_{1}, k_{2},\dots, k_{i}$.


\section{Numerical Results}
\label{numerical}

The algorithm defined in section \ref{algorithm} was implemented using 
the Java programming language. The source code can be downloaded from 
\url{http://143.107.59.106:9620/camille/beurling.tar.bz2}. This archive
provides, in fact, two slightly different implementations of the
algorithm.


\subsection{The class {\tt Beurling}}

The methods in this class compute the sequences $k_{i}$ and $b_{i}$,
storing the values in an array with fixed size.
The main method takes a natural number $n$ as input. Its output is a
file containing 3-uplas $(k_{i},\mu(k_{i}),t_{i})$, for $i$ from $3$
to $n$. Here, $t_{i}$ is the running time, in seconds, between the
computation of $\mu(k_{i-1})$ and $\mu(k_{i})$. The values of the
M\"obius function are calculated using the identity in Lemma
\ref{lemma-3} (d). The main method in this class also verifies: 
\begin{itemize}

\item
If each one of the numbers $k_{i}$ is square-free.

\item
If there exist eventually square-free numbers between $k_{i}$ and
$k_{i+1}$ that are not generated by the algorithm.

\item
The condition $k_{i+1} < 2k_{i}$, equivalent to the gap $k_{i+1} -
k_{i} <  k_{i}$.

\end{itemize}
These additional verifications are not included in the running time
$t_{i}$. \\

\begin{figure}[htb]

\begin{minipage}[b]{.46\linewidth}
\centering
\includegraphics[width=7cm]{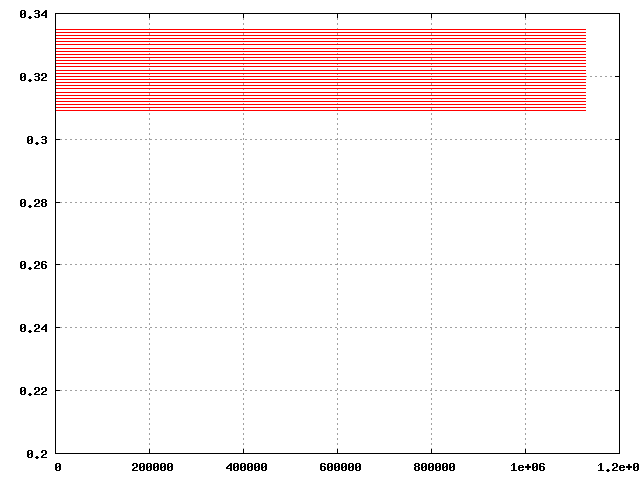}
\caption{Pairs ($k_{i}$,$t_{i}$), where $t_{i}$ is
  the running time (seconds) between the computation of
  $\mu(k_{i-1})$ and $\mu(k_{i})$, using the class {\tt Beurling}.}
\label{fig-1}
\end{minipage}
\hfill
\begin{minipage}[b]{.46\linewidth}
\centering
\includegraphics[width=7cm]{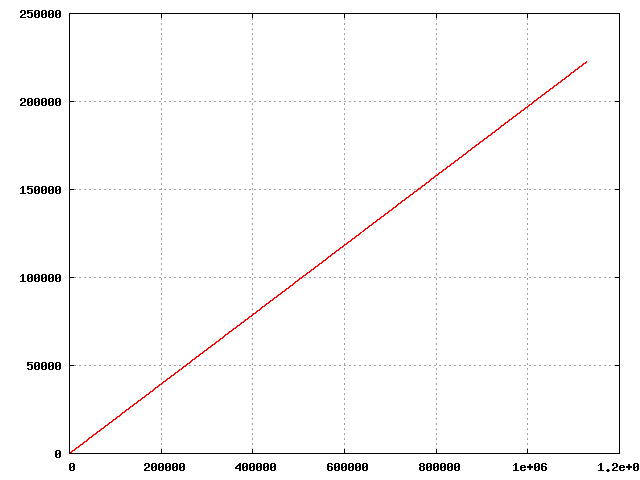}
\caption{Pairs ($k_{i}$,$T_{i}$), where $T_{i}$ is
  the \emph{total} running time (seconds) for the computation of
  $\mu(k_{i})$, using the Java class {\tt Beurling}.}
\label{fig-1b}
\vfill
\end{minipage}

\end{figure}

Running this class with $n=1\times 10^{6}$, the size of the output
file was about 10 MB. 
To estimate the running time, we generate a graphic with the
pairs $(k_{i},t_{i})$; see figure \ref{fig-1}.

This graphic seems to suggest that square-free numbers can be divided
in classes, and for each of these classes, the running time starting
from the previous iteration $t_{i}$ is constant.
We verify that all the square-free numbers in the range analyzed are
generated by the program, and satisfy the condition
$k_{i+1}<2k_{i}$, for $i \geqslant 2$. To estimate the \emph{total}
running time $T_{i}$ to compute
$\mu(k_{i})$, we use the formula $T_{i} =
\sum_{k=1}^{i}t_{k}$. The graphic with the pairs $(k_{i},T_{i})$ is
shown in figure \ref{fig-1b}. \\

The program produced also four {\it square-full}, i.e. non
square-free, numbers, shown in the table \ref{table-1}. All these
square-full numbers $k_{i}$ satisfy $b_{i}(k_{i+1}) - b_{i}(k_{i}) =
\mu(k_{i}) = 0$. However, as already stated, by the definition of the
algorithm in section \ref{algorithm}, it must be $b_{i}(k_{i+1}) -
b_{i}(k_{i}) \neq 0$ for \emph{any} generated number. Therefore, we
strongly believe that the square-full numbers are generated purely by
rounding errors.
\begin{table}[h!]
\centering
\begin{tabular}{|l|c|c|}
\hline
Number $k$ & $\mu(k)$ & Divisible by \\
\hline\hline
440375 & 0 & $5^{2}$ \\
\hline
551208 & 0 & $2^{2}$ \\
\hline
799460 & 0 & $2^{2}$ \\
\hline
979275 & 0 & $5^{2}$ \\
\hline
\end{tabular}
\caption{The four square-full numbers produced by the class {\tt
    Beurling}, with $n=1\times 10^{6}$.}
\label{table-1}
\end{table}


\subsection{The class {\tt BeurlingArrayList}}

\begin{figure}[htb]

\begin{minipage}[b]{.46\linewidth}
\centering
\includegraphics[width=7cm]{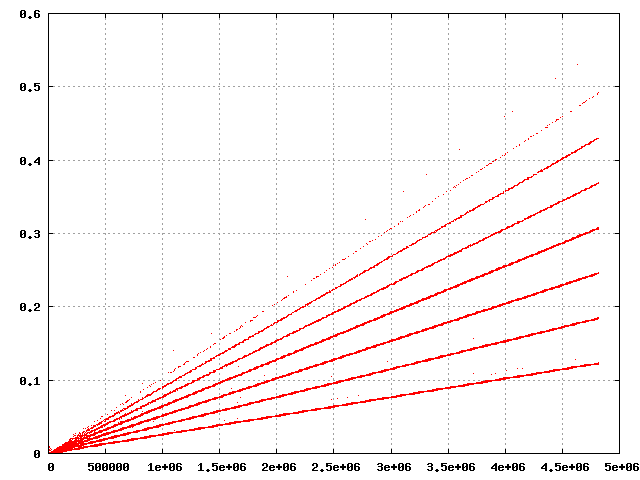}
\caption{Pairs ($k_{i}$,$t_{i}$), where $t_{i}$ is
  the running time (seconds) between the computation of
  $\mu(k_{i-1})$ and $\mu(k_{i})$, using the class {\tt
    BeurlingArrayList}.}
\label{fig-2}
\end{minipage}
\hfill
\begin{minipage}[b]{.46\linewidth}
\centering
\includegraphics[width=7cm]{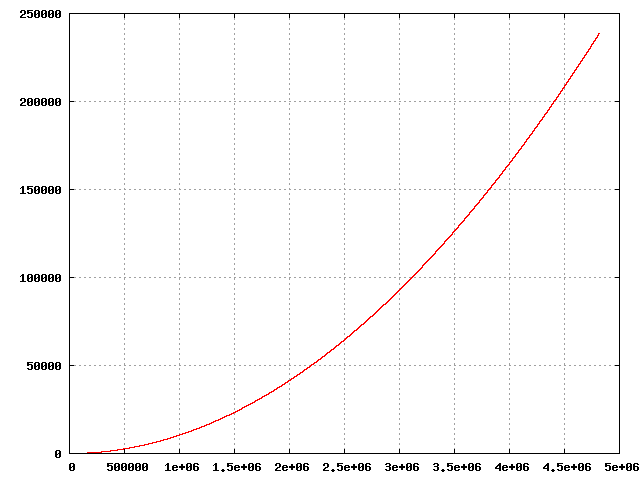}
\caption{Pairs ($k_{i}$,$T_{i}$), where $T_{i}$ is
  the \emph{total} running time (seconds) for the computation of
  $\mu(k_{i})$, using the class {\tt BeurlingArrayList}.}
\label{fig-2b}
\end{minipage}

\end{figure}

This class works as the previous one, but the computed
values are stored in an array with dynamic size. We use this class to
process $n=5\times 10^{6}$ numbers, using a relatively modest desktop
machine. The graphics with the pairs $(k_{i},t_{i})$ and
$(k_{i},T_{i})$ are shown in figures \ref{fig-2} and \ref{fig-2b},
respectively.

As in the case of figure \ref{fig-1}, the square-free numbers $k_{i}$
seems to be divided in classes. However, in this case, the running
time between two successive iterations $t_{i}$ grows linearly with
$k_{i}$. All the
square-free numbers in the range analyzed are generated by the
program, and satisfy the condition $k_{i+1}<2k_{i}$, for $i \geqslant
2$. In the range analized, twenty square-full numbers were also
generated by this class, probably by rounding errors, as explained
above.



\section{Concluding Remarks}
\label{concluding}

In section \ref{algorithm} we define an algorithm that iteratively
produces a sequence of numbers $k_{i}$ and functions $b_{i}$. 
The lemma \ref{lemma-3} states a set of properties of these sequences
suggesting that
\begin{itemize}
\item
The numbers $k_{i}$ generated by the algorithm are square-free.
\item
The set of \emph{all} the square-free numbers can be generated by the
algorithm.
\item
The value of the M\"obius function $\mu(k_{i})$ can be evaluated as
$\mu(k_{i}) = b_{i}(k_{i+1}) - b_{i}(k_{i})$.
\end{itemize}
In section \ref{main} we prove the logical equivalence of these
properties. Unfortunately, using the definition of the algorithm, we
cannot prove, neither disprove, if one of these conditions
are satisfied. Note also that in order to determine $\mu(k_{i})$,
it is necessary to generate the whole sequence $k_{1}, k_{2},\dots,
k_{i}$. \\

Numerical evidence seems to support the conjectures
quoted above. However, this evidence is limited, and certainly not
conclusive.




\vspace{1.5cm}

\noindent
Fernando Auil\\

\noindent
Escola de Artes, Ci\^encias e Humanidades\\
Universidade de S\~ao Paulo\\
Arlindo Bettio 1000\\
CEP 03828-000\\
S\~ao Paulo - SP\\
Brasil\\

\noindent
E-mail: {\tt auil@usp.br}


\end{document}